\documentclass[11pt,twoside]{article}
\usepackage{amsmath, amsthm}
\usepackage{amssymb, multicol}
\usepackage{fancyhdr, graphics, graphicx}

\title{\bf \sc{Representation of Period Doubling by Digraphs and Characteristic Polynomials}}
\vspace{.2in}
\author{
{\small{\textsc{Richell O. Celeste}}}\\
{\small{Institute of Mathematics}}\\
{\small{University of the Philippines}}\\
{\small{Diliman Quezon City}}\\
\small e-mail: \texttt{ching@math.upd.edu.ph} \\
\\
{\small{\textsc{Yoshifumi Takenouchi}}}\\
{\small{Institute of Mathematics}}\\
{\small{University of the Philippines}}\\
{\small{Diliman Quezon City}}\\
\small e-mail: \texttt{jytakenovich@math.upd.edu.ph}
}
\date{ }

\newtheorem{thm}{Theorem}[section]
\newtheorem{lem}{Lemma}[section]

\newtheorem{remark}{Remark}[section]
\newtheorem{Def}{Definition}[section]
\newtheorem{Cor}{Corollary}[section]

\begin{document}

\maketitle

\begin{abstract}

   A general procedure which defines a partial ordering of cyclic permutations induced by continuous maps
is known for constructing immediate successors to a cycle. We expound on this procedure in terms of labelled digraphs and characteristic polynomials then apply this study to period doubling, the most common route to chaos for a nonlinear dynamical system.\\

\end{abstract}

\section{Introduction}
%\subparagraph{earlier literatures}
$\ \ \ \ $
%\subparagraph{1.Li and Yorke}
Li and Yorke \cite{LY} showed that
the presence of a periodic orbit with least period $3$ forces
the presence of periodic orbit with all possible least periods.
%\subparagraph{2.Sharkovsky}
However before \cite{LY}, Sharkovsky \cite{BH,Drazin} had already obtained the results
which include the statement above as a corollary.
%\subparagraph{3.Baldwin}
Baldwin \cite{B} extended this idea to
get a partial order on the set of cycles
in which the type of the periodic orbit is taken into account.
%\subparagraph{4.Bernhardt}
Bernhardt \cite{Bernhardt1} showed how to construct immediate successors to a given cycle.

%\subparagraph{standing of this paper among them}
First we reexamine \cite{Bernhardt1} to
obtain graph theoretical explanation of
the relation between a cycle and its
immediate successors.
With this view point, we then obtain
the representation of the period-doubling bifurcation phenomena
by the digraphs consisting of signed vertices and directed edges and
by the characteristic polynomials of the adjacency matrices.
Obtaining a discrete expression of the period-doubling bifurcation phenomena
implies that we may define the period-doubling bifurcation phenomena on the graph itself
without referring to the original one-dimensional continuous maps.

%\subparagraph{construction of this paper}
In section 2 we present basic definitions and concepts
on the forcing relation on cycles in dimension one and
on graph theory.
In section 3 we obtain the relation between the digraph and the characteristic polynomial of any cycle and that of its immediate successor.
In section 4 we describe the period-doubling phenomena of the family of the logistic curves with their digraphs and characteristic polynomials.

\section{Preliminaries}

%\subparagraph{period, orbit}
\begin{Def}\rm{\bf }
Let $f:I\to I$ be a continuous map of a compact interval $I$ to itself.
We define $f^0(x)=x$ and for $n\in\mathbb{N} \setminus\{0\}, f^n(x)=f(f^{n-1}(x))$.
If there exists $k\in\mathbb{N}$ such that $f^k(x)=x$,
then we say $x$ is \textbf{periodic} for $f$ and
$x$ has (least) \textbf{period} $s$,
where $s$ is smallest element of $\mathbb{N}$ such that $f^s(x)=x$.
For $x\in I$, the \textbf{orbit} of $x$ is the set $\{f^n(x)|n\ge0\}$.
If $x$ is periodic with least period $s$, then the orbit of $x$ is the finite set
${\cal O}(x)=\{x,f(x),f^2(x),\ldots,f^{s-1}(x)\}$.
\end{Def}

%\subparagraph{cycle, type}
\begin{Def}\rm{\bf }
A \textbf{cycle} (cyclic permutation) of length $n$ is
a bijection $\theta:\{1,2,\ldots,n\}$
$\to\{1,2,\ldots,n\}$
such that $\theta^k(1)\neq1$ for $1\le k< n$.
We denote a cycle $\theta$ by $\theta=(k_1,k_2,\ldots,k_m)$,
where $\theta(k_i)=k_{i+1}$ and $\theta(k_m)=k_1$.
We assume, without loss of generality, that $k_1=1$.
Write the elements of a periodic orbit ${\cal O}$ in increasing order:
$x_1<x_2<\cdots<x_s$.
We say ${\cal O}$ has \textbf{orbit type} $\theta$ if
$\theta$ is a cycle of length $s$ and for each
$i\in\{1,2,\ldots,s\},\ f(x_i)=x_{\theta(i)}$.
\end{Def}
$(S_n, \circ)$ denotes the group of permutations on $n$ objects,
$C_n$ denotes the subset of $S_n$ consisting of cycles (cyclic permutations) of length $n$.
We also denote $\bigcup_{n\ge1}C_n$ by $C$.

%\subparagraph{forcing relation}
\begin{Def}\rm{\bf }
The \textbf{forcing relation on cycles} is defined as follows:
$\beta\in C$ forces $\alpha\in C$
if and only if
every continuous map of the interval that
has a periodic orbit of type $\beta$
has a periodic orbit of type $\alpha$.
\end{Def}
Baldwin \cite{B} showed that
this relation induces a partial order on the set of cycles and
provided an algorithm to decide
when one cycle forces another.

%\subparagraph{immediate successor}
\begin{Def}\rm{\cite{Bernhardt1}}
We say $\beta\in C$ is an \textbf{immediate successor} to $\alpha\in C$, or
$\alpha$ is an \textbf{immediate predecessor} to $\beta$,
if $\beta\neq\alpha$  and  the set of all cycles forced by $\beta$
is same as the set of all cycles forced by $\alpha$ and $\alpha$ itself.
\end{Def}

%\subparagraph{linear map}
\begin{Def}\rm{\bf }
For $\theta\in S_n$,
the \textbf{canonical $\theta$-linear map}, or the \textbf{connect-the-dots map of $\theta$}
is defined by $L_{\theta}:[1,n]\to[1,n]$ where
$L_{\theta}=\theta$ on $\{1,2,\ldots,n\}$ and
$L_{\theta}$ is linear on $I_i:=[i,i+1]$ for $i\in\{1,2,\ldots,n-1\}$.
%Here we note that the graph of $L_{\theta}$ consists of at most $n-1$ linear segments, each having a slope $m\in\mathbb{Z}$ satisfying $|m|\ge1$.
\end{Def}
%\subparagraph{well known}
It is known that
$\beta\in C$ forces $\alpha\in C$
if and only if
$L_\beta$ has a periodic orbit of cycle type $\alpha$
\cite{ALM}.

%\subparagraph{modality}
\begin{Def}
For $\theta\in C_n,\ n\ge3$
we say $\theta$ is \textbf{$+k$-modal cycle} or \textbf{$-k$-modal cycle} if
the number of all the local extrema of $L_{\theta}$ is equal to $k$
and $\theta(1)<\theta(2)$ or $\theta(1)>\theta(2)$, respectively.
In particular $+1$-modal cycle is called \textbf{unimodal cycle}.
\end{Def}
The forcing relation defined by Baldwin induces
a total order on the set of unimodal cycles \cite{CandE}.

%\subparagraph{digraph with signed vertices and directed edges}
\begin{Def}\rm{\bf }
From $L_{\theta}$,
we define the \textbf{labelled digraph (directed graph) of $\theta$}, denoted by $G(\theta)$, as follows:
$G(\theta)=({\cal V,E})$ consists of
the set of $n-1$ signed vertices
${\cal V}=\{v_{1}^{\pm},v_{2}^{\pm},\ldots,v_{n-1}^{\pm}\}$
where we take $v_{i}^{+}$ if $\theta(i)<\theta(i+1)$ or $v_{i}^{-}$ otherwise,
and the set of directed edges ${\cal E}=\{e_{ij}\}$ where $e_{ij}$
a directed edge from $v_{i}$ to $v_{j}$ exists if $I_j\subset L_{\theta}(I_i)$.
\end{Def}
Conversely, from $G(\theta)$ we can obtain $L_{\theta}$.
In this sense, we say they are equivalent.
Indeed Baldwin's algorithm \cite{B}
which decides whether one cycle $\theta$ forces another or not,
is derived from $G(\theta)$.
%\subparagraph{adjacency matrix, induced matrix}
\begin{Def}\rm{\bf }
The \textbf{induced matrix of $\theta$} or the \textbf{adjacency matrix of $G(\theta)$} \rm{\cite{graph}},
denoted by $M(\theta)=\Big[a_{i,j}\Big]$, is the $(n-1)\times(n-1)$
matrix with $(i,j)$-th entry given by
$$
a_{i,j}=
\left\{
\begin{array}{cl}
1,&\mbox{if $I_j\subset L_{\theta}(I_i)$,}\\
0,&\mbox{otherwise,}\\
\end{array}
\right.
$$
where $i, j \in[1, n-1]$.
\end{Def}

%\subparagraph{induced subgraph}
\begin{Def}\rm{\cite{graph}}
A graph $H=({\cal W},{\cal F})$ is said to be a \textbf{subgraph} of
the graph $G=({\cal V},{\cal E})$
if ${\cal W}\subset{\cal V}$ and ${\cal F}\subset{\cal E}$.
In particular if ${\cal F}$ consists of all the edges from ${\cal E}$
which connect the vertices from ${\cal W}$,
then $H$ is called an \textbf{induced subgraph}.
\end{Def}

%\subparagraph{cycle}
\begin{Def}\rm{\cite{graph}}
A \textbf{cycle of length $n$} in the digraph $G$
is a subgraph with the vertex set $\{v_1,\ldots,v_n\}$
having directed edges from $v_i$ to $v_{i+1}$ where $i=1,\ldots,n-1$ and
a directed edge from $v_n$ to $v_{1}$. In particular if $v_1<v_2<\cdots<v_n$ then this cycle induces the cyclic permutation $\theta=(12\cdots n)$.
\end{Def}

\bigskip

\section{Theorem}
%\subparagraph{theta^*}
\begin{Def}\rm{\cite{Bernhardt1}}
Let $\theta\in C_n$,
then $\theta^*\in S_{2n}$ is defined by
$$
\theta^*(2k):=2\theta(k),\ \theta^*(2k-1):=2\theta(k)-1,
$$
where $k\in\{1,2,\ldots,n\}$.
\end{Def}
%\subparagraph{\eta}
\begin{Def}\rm{\cite{Bernhardt1}}
Let $\rho_s$ denote the transposition $(2s-1, 2s)$, where $s\in\mathbb{N}$.
Then for $\theta\in C_n$ define $\eta$ as follows:
$$
\eta:=\theta^*\circ\underbrace{
\rho_{i_1}\circ\rho_{i_2}\circ\cdots\circ\rho_{i_{m}}
}_{m},
$$
where $1\le i_1<i_2<\cdots<i_m\le n$.
\end{Def}
If $m$ is odd, then $\eta\in C_{2n}$ and it is an immediate successor to $\theta$.
If $m$ is even, then $\eta\in S_{2n}\setminus C_{2n}$.
There are $2^{n-1}$ immediate successors $\eta\in C_{2n}$ to $\theta$,
and $2^{n-1}$ non-cyclic permutations $\eta\in S_{2n}$. (See \rm{\cite{Bernhardt0}} and \rm{\cite{Bernhardt1}})
\\

%\subsubsection{digraph}
For $\theta$ and $\eta$ above,
denote their labelled digraphs by $G(\theta)$ and $G(\eta)$,
respectively.
Then $G(\theta)=({\cal V,E})$ consists of
the set of $n-1$ signed vertices
${\cal V}=\{v_{1}^{\pm},v_{2}^{\pm},\ldots,v_{n-1}^{\pm}\}$,
where we take $v_{i}^{+}$ if $\theta(i)<\theta(i+1)$ or $v_{i}^{-}$ otherwise,
and
the set of directed edges ${\cal E}=\{e_{ij}\}$ where $e_{ij}$
a directed edge from $v_{i}$ to $v_{j}$ exists if $I_j\subset L_{\theta}(I_i)$.
$G(\eta)=({\cal W},{\cal F})$ consists of
the set of $2n-1$ signed vertices
${\cal W}=\{w_{1}^{\pm},w_{2}^{\pm},\ldots,w_{2n-1}^{\pm}\}$,
where we take $w_{i}^{+}$ if $\eta(i)<\eta(i+1)$ or $w_{i}^{-}$ otherwise,
and
the set of directed edges ${\cal F}=\{f_{ij}\}$ where $f_{ij}$
a directed edge from $w_{i}$ to $w_{j}$ exists if $I_j\subset L_{\eta}(I_i)$.
Let us divide ${\cal W}$ into two subsets
${\cal W}_1=\{w_{2}^{\pm},\ldots,w_{2n-2}^{\pm}\}$
and
${\cal W}_2=\{w_{1}^{\pm},w_{3}^{\pm},\ldots,w_{2n-1}^{\pm}\}$.
\\
\\
Then we have the following observations:
\\
\\
%\subparagraph{i}
$(i)$ $G(\eta)|_{{\cal W}_1}$, the induced subgraph of $G(\eta)=(\cal{W},\cal{F})$
determined on ${\cal W}_1$ is isomorphic to
the digraph $G(\theta)=({\cal V,E})$.
Indeed, the isomorphism $\varphi:{\cal V}\to{\cal W}_1$
is given by
$$
\varphi(v_{i}^{\pm})=w_{2i}^{\pm}
$$
where $i\in\{1,2,\ldots,n-1\}$.
Observe that
the sign of vertex $v_{i}$ is always consistent with that of vertex $w_{2i}$.
\\
\\
%\subparagraph{ii}
$(ii)$ In $G(\eta)=(\cal{W},\cal{F})$
there does not exist any directed edge from ${\cal W}_2$ to ${\cal W}_1$.
\\
\\
%\subparagraph{iii}
$(iii)$ $G(\eta)|_{{\cal W}_2}$, the induced subgraph of $G(\eta)=(\cal{W},\cal{F})$
determined on ${\cal W}_2$, consists of an $n$-cycle induced by $\theta$.
Indeed, this $n$-cycle is given by
\begin{center}
\unitlength=0.8mm
\begin{picture}(100,25)
\put(  0,10){\circle*{2}}
\put( 20,10){\circle*{2}}
\put( 40,10){\circle*{2}}
\put( 60,10){}
\put( 80,10){}
\put(100,10){\circle*{2}}
\put(-1,5){$w_{1}$}
\put(19,5){$w_{2\theta(1)-1}$}
\put(39,5){$w_{2\theta^2(1)-1}$}
\put(59,5){}
\put(79,5){}
\put(99,5){$w_{2\theta^{n-1}(1)-1}$}
\put( 2,10){\vector(1,0){16}}
\put(22,10){\vector(1,0){16}}
\put(42,10){\vector(1,0){16}}
\put(62,10){$\ldots\ldots\ldots$}
\put(82,10){\vector(1,0){16}}
\put(100,20){\line(-1,0){100}}
\put(100,20){\line(0,-1){8}}
\put(  0,20){\vector(0,-1){8}}
\end{picture}
\end{center}
and the type of this $n$-cycle is also $\theta$.
\\

Together, $(i)$, $(ii)$, and $(iii)$ provide a graph theoretical explanation of the reason why
$\eta$ is an immediate successor to $\theta$
(see {\textsc{Lemma 2.4}} in \rm{\cite{Bernhardt1}}).
As a consequence, we have the following theorem.

%\subsubsection{Main Thm}
\begin{thm}
Denote the characteristic polynomials of the adjacency matrices of
$G(\theta),G(\eta),G(\eta)|_{{\cal W}_1}$ and $G(\eta)|_{{\cal W}_2}$ by
$P_{\theta}(\lambda),P_{\eta}(\lambda),P_{{\cal W}_1}(\lambda)$ and $P_{{\cal W}_2}(\lambda)$, respectively.
Then $P_{\eta}(\lambda)$ is given by the product of $P_{{\cal W}_1}(\lambda)$ and $P_{{\cal W}_2}(\lambda)$,
namely
$$
P_{\eta}(\lambda)=P_{\theta}(\lambda)\cdot(\lambda^n-1).
$$
\end{thm}
As a corollary of this theorem, we have the following.
%\subparagraph{Cor}
\begin{Cor}
The spectrum of $G(\eta)$ consists of
the spectrum of $G(\theta)$ and the $n$ distinct roots of unity.
\end{Cor}

Assume that the adjacency matrix of $G(\theta)$ is given by
$$
M(\theta)=\Bigg[\ a_{i,j}\ \Bigg]_{(n-1)\times(n-1)},
$$
where $1\le i,j\le n-1$,
and the adjacency matrix of $G(\eta)$ is given by
$$
M(\eta)=\Bigg[\ b_{i,j}\ \Bigg]_{(2n-1)\times(2n-1)},
$$
where $1\le i,j\le 2n-1$.

Then we have the following Lemma 3.1, and 3.2.

%\subsubsection{lemma 1}
\begin{lem} $G(\eta)|_{{\cal W}_1}$ is isomorphic to
$G(\theta)$ by $\varphi$, namely
$$
b_{2i,2j}=a_{i,j},
$$
where $1\le i,j\le n-1$.
\end{lem}
%\subparagraph{Proof}
\textbf{Proof.}
For any $\theta\in C_n$, $\theta^*\in S_{2n}$ is defined by
$$
\theta^*:=
\left(
  \begin{array}{cccccc}
\cdots&        2i-1&        2i&          2i+1&        2i+2&\cdots\\
\cdots&2\theta(i)-1&2\theta(i)&2\theta(i+1)-1&2\theta(i+1)&\cdots\\
  \end{array}
\right).
$$
If $\eta:=\theta^*\circ\cdots\circ\rho_{i}\circ\cdots$, then
$$
\eta=
\left(
  \begin{array}{cccccc}
\cdots&      2i-1&          2i&          2i+1&        2i+2&\cdots\\
\cdots&2\theta(i)&2\theta(i)-1&2\theta(i+1)-1&2\theta(i+1)&\cdots\\
  \end{array}
\right),
$$
and
$$
L_{\eta}([2i,2i+1])=
\left\{
\begin{array}{ccc}
{\displaystyle [2\theta (i)-1,2\theta (i+1)-1]}&\mbox{if}&\theta (i)<\theta (i+1)\\
{\displaystyle [2\theta (i+1)-1,2\theta (i)-1]}&\mbox{if}&\theta (i)>\theta (i+1)\\
\end{array}
\right..
$$
If $\eta:=\theta^*\circ\cdots\circ\rho_{i+1}\circ\cdots$, then
$$
\eta=
\left(
  \begin{array}{cccccc}
\cdots&        2i-1&        2i&        2i+1&          2i+2&\cdots\\
\cdots&2\theta(i)-1&2\theta(i)&2\theta(i+1)&2\theta(i+1)-1&\cdots\\
  \end{array}
\right),
$$
and
$$
L_{\eta}([2i,2i+1])=
\left\{
\begin{array}{ccc}
{\displaystyle [2\theta (i),2\theta (i+1)]}&\mbox{if}&\theta (i)<\theta (i+1)\\
{\displaystyle [2\theta (i+1),2\theta (i)]}&\mbox{if}&\theta (i)>\theta (i+1)\\
\end{array}
\right..
$$
If $\eta:=\theta^*\circ\cdots\circ\rho_{i}\circ\rho_{i+1}\circ\cdots$, then
$$
\eta=
\left(
  \begin{array}{cccccc}
\cdots&        2i-1&        2i&        2i+1&          2i+2&\cdots\\
\cdots&2\theta(i)&2\theta(i)-1&2\theta(i+1)&2\theta(i+1)-1&\cdots\\
  \end{array}
\right),
$$
and
$$
L_{\eta}([2i,2i+1])=
\left\{
\begin{array}{ccc}
{\displaystyle [2\theta (i)-1,2\theta (i+1)]}&\mbox{if}&\theta (i)<\theta (i+1)\\
{\displaystyle [2\theta (i+1),2\theta (i)-1]}&\mbox{if}&\theta (i)>\theta (i+1)\\
\end{array}
\right..
$$
If $\eta:=\theta^*\circ\rho_{i_1}\circ\rho_{i_2}\circ\cdots\circ\rho_{i_{m}}$ does not include
$\rho_{i}$ nor $\rho_{i+1}$, then
$$
\eta=
\left(
  \begin{array}{cccccc}
\cdots&        2i-1&        2i&          2i+1&        2i+2&\cdots\\
\cdots&2\theta(i)-1&2\theta(i)&2\theta(i+1)-1&2\theta(i+1)&\cdots\\
  \end{array}
\right),
$$
and
$$
L_{\eta}([2i,2i+1])=
\left\{
\begin{array}{ccc}
{\displaystyle [2\theta (i),2\theta (i+1)-1]}&\mbox{if}&\theta (i)<\theta (i+1)\\
{\displaystyle [2\theta (i+1)-1,2\theta (i)]}&\mbox{if}&\theta (i)>\theta (i+1)\\
\end{array}
\right..
$$
Under the assumption $\theta(i)<\theta(i+1)$, we have
$$
\begin{array}{rl}
               &a_{i,j}=1\\
\Leftrightarrow&I_j\subset L_{\theta}(I_i)\\
\Leftrightarrow&[j,j+1]\subset[\theta(i),\theta(i+1)]\\
\Leftrightarrow&\theta(i) \le j\mbox{ and }j+1 \le \theta(i+1)\\
\Leftrightarrow&2\theta(i) \le 2j\mbox{ and }2j+1 \le 2\theta(i+1)-1\\
\Leftrightarrow&[2j,2j+1]\subset[2\theta(i),2\theta(i+1)-1].\\
\end{array}
$$
Here we note that
$$
[2\theta(i),2\theta(i+1)-1]\subset[2\theta(i)-1,2\theta(i+1)-1],[2\theta(i),2\theta(i+1)]\subset[2\theta(i)-1,2\theta(i+1)],
$$
and
$$
[{\displaystyle \mathop{2j}_{\mbox{\tiny{even}}}},
{\displaystyle \mathop{2j+1}_{\mbox{\tiny{odd}}}}]\neq
[{\displaystyle \mathop{2\theta(i)-1}_{\mbox{\tiny{odd}}}},
{\displaystyle \mathop{2\theta(i)}_{\mbox{\tiny{even}}}}],\
[{\displaystyle \mathop{2\theta(i+1)-1}_{\mbox{\tiny{odd}}}},
{\displaystyle \mathop{2\theta(i+1)}_{\mbox{\tiny{even}}}}],
$$
thus,
$$
\begin{array}{rl}
               &[2j,2j+1]\subset[2\theta(i),2\theta(i+1)-1]\\
\Leftrightarrow&[2j,2j+1]\subset[\eta(2i),\eta(2i+1)]\\
\Leftrightarrow&I_{2j}\subset L_{\eta}(I_{2i})\\
\Leftrightarrow&b_{2i,2j}=1.\\
\end{array}
$$
Similarly,
under the assumption $\theta(i)>\theta(i+1)$, we can also show that
$$
a_{i,j}=1\Leftrightarrow b_{2i,2j}=1.
$$
This concludes the proof of the lemma.
$\hspace*{\fill}\Box$

%\subsubsection{lemma 2}
\begin{lem} In $G(\eta)=(\cal{W},\cal{F})$, there does not exist any directed edge from ${\cal W}_2$ to ${\cal W}_1$, that is
$$
b_{2i-1,2j}=0,
$$
where $1\le i\le n,\ 1\le j\le n-1$.
Moreover, $G(\eta)|_{{\cal W}_2}$ consists of an $n$-cycle induced by $\theta$,
namely
$$
\left\{
\begin{array}{cl}
b_{2i-1,2j-1}=1&\mbox{if $j=\theta(i)$,}\\
b_{2i-1,2j-1}=0&\mbox{otherwise,}\\
\end{array}
\right.
$$
where $1\le i,j\le n$.

In other words
$$
\left[
\begin{array}{cccc}
b_{12}     &b_{14}     &\cdots&b_{1(2n-2)}       \\
b_{32}     &b_{34}     &\cdots&b_{3(2n-2)}       \\
\cdots     &\cdots     &\cdots&\cdots            \\
\cdots     &\cdots     &\cdots&\cdots            \\
b_{(2n-1)2}&b_{(2n-1)4}&\cdots&b_{(2n-1)(2n-2)}  \\
\end{array}\right]=O_{n\times(n-1)},
$$
and
$$
\left[
\begin{array}{cccc}
b_{11}     &b_{13}     &\cdots&b_{1(2n-1)}     \\
b_{31}     &b_{33}     &\cdots&b_{3(2n-1)}     \\
\cdots     &\cdots     &\cdots&\cdots          \\
\cdots     &\cdots     &\cdots&\cdots          \\
b_{(2n-1)1}&b_{(2n-1)3}&\cdots&b_{(2n-1)(2n-1)}\\
\end{array}
\right]
=
\left[
\begin{array}{cccc}
c_{11}&c_{12}&\cdots&c_{1n}\\
c_{21}&c_{22}&\cdots&c_{2n}\\
\cdots&\cdots&\cdots&\cdots\\
\cdots&\cdots&\cdots&\cdots\\
c_{n1}&c_{n2}&\cdots&c_{nn}\\
\end{array}
\right]
$$
$$
\Longleftrightarrow\ b_{2i-1,2j-1}=c_{ij},
$$
where
$$
\left\{
\begin{array}{cl}
c_{ij}=1&\mbox{if $j=\theta(i)$,}\\
c_{ij}=0&\mbox{otherwise.}\\
\end{array}
\right.
$$
\end{lem}
%\subparagraph{Proof}
\textbf{Proof.}
$$
\begin{array}{rl}
               &j=\theta(i)\\
\Leftrightarrow&\theta^*(2i-1)=2j-1\mbox{ and }\theta^*(2i)=2j\\
\Leftrightarrow&[2j-1,2j]=[\theta^*(2i-1),\theta^*(2i)].\\
\end{array}
$$
For any $\theta\in C_n$, $\theta^*\in S_{2n}$ is defined by
$$
\theta^*:=
\left(
  \begin{array}{cccc}
\cdots&        2i-1&        2i&\cdots\\
\cdots&\theta^*(2i-1)&\theta^*(2i)&\cdots\\
  \end{array}
\right).
$$
If $\eta:=\theta^*\circ\cdots\circ\rho_{i}\circ\cdots$, then
$$
\eta=
\left(
  \begin{array}{cccc}
\cdots&        2i-1&        2i&\cdots\\
\cdots&\theta^*(2i)&\theta^*(2i-1)&\cdots\\
  \end{array}
\right),
$$
and thus $[2j-1,2j]=[\theta^*(2i-1),\theta^*(2i)]=[\eta(2i),\eta(2i-1)]$.

If $\eta:=\theta^*\circ\rho_{i_1}\circ\rho_{i_2}\circ\cdots\circ\rho_{i_{m}}$ does not include $\rho_{i}$, then
$$
\eta=
\left(
  \begin{array}{cccc}
\cdots&        2i-1&        2i&\cdots\\
\cdots&\theta^*(2i-1)&\theta^*(2i)&\cdots\\
  \end{array}
\right),
$$
and thus $[2j-1,2j]=[\theta^*(2i-1),\theta^*(2i)]=[\eta(2i-1),\eta(2i)]$.

In either event, we have
$$
\begin{array}{rl}
               &[2j-1,2j]=[\theta^*(2i-1),\theta^*(2i)]\\
\Leftrightarrow&I_{2j-1}= L_{\eta}(I_{2i-1})\\
\Leftrightarrow&\left\{
\begin{array}{ll}
b_{2i-1,2j-1}=1&,\\
b_{2i-1,k}=0&\mbox{if $k\neq2j-1$,}\\
\end{array}
\right.
\\
\end{array}
$$
ending the proof of the lemma.
$\hspace*{\fill}\Box$

The following statement is obvious:
%\subsubsection{lemma 3}
\begin{lem}
$$
\left|
\begin{array}{cc}
A_{l\times l}&O_{l\times m}\\
C_{m\times l}&B_{m\times m}\\
\end{array}
\right|=|A_{l\times l}||B_{m\times m}|.
$$
\end{lem}

We are now ready to prove the theorem.
\\
\\
%\subsubsection{Proof of the thm}
\textbf{Proof of Theorem 3.1.}
The characteristic polynomial of the adjacency matrix of $G(\theta)$ is given by
$$
P_{\theta}(\lambda)=|\lambda I-M(\theta)|=\left|
\begin{array}{cccc}
\lambda-a_{11}&       -a_{12}&\cdots&           -a_{1(n-1)}\\
       -a_{21}&\lambda-a_{22}&\cdots&           -a_{2(n-1)}\\
        \cdots&        \cdots&\cdots&                \cdots\\
        \cdots&        \cdots&\cdots&                \cdots\\
        \cdots&        \cdots&\cdots&                \cdots\\
   -a_{(n-1)1}&   -a_{(n-1)2}&\cdots&\lambda-a_{(n-1)(n-1)}\\
\end{array}
\right|.
$$
Let
$$
\sigma_{2n-1}:=\left(
\begin{array}{cccccccccc}
1&2&3&\cdots&   n&n+1&n+2&n+3&\cdots&2n-1\\
1&3&5&\cdots&2n-1&  2&  4&  6&\cdots&2n-2\\
\end{array}
\right).
$$
Using Lemmas 3.1, 3.2, and 3.3,  the characteristic polynomial  $P_{\eta}(\lambda)$ of $\eta$ is thus
$$
\begin{array}{rcl}
&=&\left|
\begin{array}{ccccccc}
\lambda-b_{11}&-b_{12}       &\cdots&\cdots&-b_{1(2n-2)}      &-b_{1(2n-1)}\\
-b_{21}       &\lambda-b_{22}&\cdots&\cdots&-b_{2(2n-2)}      &-b_{2(2n-1)}\\
\cdots        &\cdots        &\cdots&\cdots&\cdots            &\cdots\\
\cdots        &\cdots        &\cdots&\cdots&\cdots            &\cdots\\
-b_{(2n-2)1}  &-b_{(2n-2)2}  &\cdots&\cdots&\lambda-b_{(2n-2)(2n-2)}&-b_{(2n-2)(2n-1)}\\
-b_{(2n-1)1}  &-b_{(2n-1)2}  &\cdots&\cdots&-b_{(2n-1)(2n-2)}     &\lambda-b_{(2n-1)(2n-1)}\\
\end{array}
\right|\\
\\
&=&\mbox{sgn}(\sigma_{2n-1})\left|
\begin{array}{ccccccc}
\lambda-b_{11}&-b_{12}       &\cdots&\cdots&-b_{1(2n-2)}      &-b_{1(2n-1)}\\
\cdots        &\cdots        &\cdots&\cdots&\cdots            &\cdots\\
-b_{(2n-2)1}  &-b_{(2n-2)2}  &\cdots&\cdots&\lambda-b_{(2n-2)(2n-2)}&-b_{(2n-2)(2n-1)}\\
-b_{21}       &\lambda-b_{22}&\cdots&\cdots&-b_{2(2n-2)}      &-b_{2(2n-1)}\\
\cdots        &\cdots        &\cdots&\cdots&\cdots            &\cdots\\
-b_{(2n-1)1}  &-b_{(2n-1)2}  &\cdots&\cdots&-b_{(2n-1)(2n-2)} &\lambda-b_{(2n-1)(2n-1)}\\
\end{array}
\right|\\
\\
&=&(\mbox{sgn}(\sigma_{2n-1}))^2\left|
\begin{array}{ccccccc}
\lambda-b_{11}&\cdots&-b_{1(2n-1)}      &-b_{12}     &\cdots&-b_{1(2n-2)}      \\
\cdots        &\cdots&\cdots            &\cdots      &\cdots&\cdots            \\
-b_{(2n-1)1}  &\cdots&\lambda-b_{(2n-1)(2n-1)}&-b_{(2n-1)2}&\cdots&-b_{(2n-1)(2n-2)} \\
-b_{21}       &\cdots&-b_{2(2n-1)}      &\lambda-b_{22}&\cdots&-b_{2(2n-2)}      \\
\cdots        &\cdots&\cdots            &\cdots      &\cdots&\cdots            \\
-b_{(2n-2)1}  &\cdots&-b_{(2n-2)(2n-1)} &-b_{(2n-2)2}&\cdots&\lambda-b_{(2n-2)(2n-2)}\\
\end{array}
\right|\\
\\
&=&\left|
\begin{array}{cccc}
\lambda-b_{11}&\cdots&-b_{1(2n-1)}      \\
\cdots        &\cdots&\cdots            \\
-b_{(2n-1)1}  &\cdots&\lambda-b_{(2n-1)(2n-1)}\\
\end{array}
\right|
\left|
\begin{array}{cccc}
\lambda-b_{22}&\cdots&-b_{2(2n-2)}      \\
\cdots        &\cdots&\cdots            \\
-b_{(2n-2)2}  &\cdots&\lambda-b_{(2n-2)(2n-2)}\\
\end{array}
\right|\\
\\
&=&(\lambda^n+\mbox{sgn}(\theta)(-1)^n)P_{\theta}(\lambda)\\
\end{array}
$$
since
$
(\mbox{sgn}(\sigma_{2n-1}))^2=1,
$
\\
$$
\left[
\begin{array}{ccc}
-b_{12}     &\cdots&-b_{1(2n-2)}       \\
\cdots     &\cdots&\cdots            \\
-b_{(2n-1)2}&\cdots&-b_{(2n-1)(2n-2)}  \\
\end{array}
\right]=O_{n\times(n-1)},
$$
\\
$$
\left|
\begin{array}{cccc}
\lambda-b_{22}&\cdots&-b_{2(2n-2)}      \\
\cdots        &\cdots&\cdots            \\
-b_{(2n-2)2}  &\cdots&\lambda-b_{(2n-2)(2n-2)}\\
\end{array}
\right|
=
\left|
\begin{array}{cccc}
\lambda-a_{11}&\cdots&-a_{1(n-1)}      \\
\cdots        &\cdots&\cdots          \\
-a_{(n-1)1}   &\cdots&\lambda-a_{(n-1)(n-1)}\\
\end{array}
\right|
=P_{\theta}(\lambda),
$$
and
$$
\begin{array}{rl}
&\left|
\begin{array}{cccc}
\lambda-b_{11}&\cdots&-b_{1(2n-1)}      \\
\cdots        &\cdots&\cdots            \\
-b_{(2n-1)1}  &\cdots&\lambda-b_{(2n-1)(2n-1)}\\
\end{array}
\right|\\
\\
=&\left|
\begin{array}{cccc}
\lambda-c_{11}&\cdots&-c_{1n}\\
\cdots        &\cdots&\cdots\\
-c_{n1}       &\cdots&\lambda-c_{nn}\\
\end{array}
\right|\mbox{ where }\left\{
\begin{array}{cl}
c_{ij}=1&\mbox{if $j=\theta(i)$,}\\
c_{ij}=0&\mbox{otherwise.}\\
\end{array}
\right.\\
\\
=&\mbox{sgn}(\iota)\cdot(\lambda-0)^n+\mbox{sgn}(\theta)\cdot(-1)^n\\
=&\lambda^n+\mbox{sgn}(\theta)(-1)^n.\\
\end{array}
$$
Here we note that $\mbox{sgn}(\theta)(-1)^n$ is always equal to $-1$ since
$$
\mbox{sgn}(\theta)(-1)^n=
\left\{
\begin{array}{cl}
(+1)\cdot(-1)&\mbox{if $n$ is odd,}\\
(-1)\cdot(+1)&\mbox{if $n$ is even.}\\
\end{array}
\right.
$$
This completes the proof of the theorem.
$\hspace*{\fill}\Box$

\section{Immediate Successors and Period Doubling}
\subsection{The logistic map}
\hspace{.45in} The logistic map \cite{Fractals and Chaos, Drazin} is written as
\begin{equation}\label{l1}
x_{n+1}=ax_n(1-x_n)
\end{equation}
where the current value $x_n$, is mapped onto the next value $x_{n+1}$. Here, we restrict our considerations to $0\leq x_n\leq1$ and $0\leq a\leq 4$. The corresponding function to the logistic map, called the logistic function or the logistic curve, is given by
 \begin{equation}
f_{a}(x)=ax(1-x)
\end{equation}
By repeatedly iterating the logistic map forward through time, we may observe different behaviors of the iterated solutions. The sequence of iterated
solutions to the map is called an orbit. The behavior of orbits, which originate from typically chosen initial condition used for the iterations, depends on the control parameter $a$. The logistic curve is parabolic like the quadratic function with $f_a(0) = f_a(1) = 0$ and a maximum at $x=1/2$. Being a smooth curve with only one critical point, that is, a local maximum,  it is classified as a unimodal (single humped) map.

\begin{center}
  \includegraphics[width=4in,height=2in]{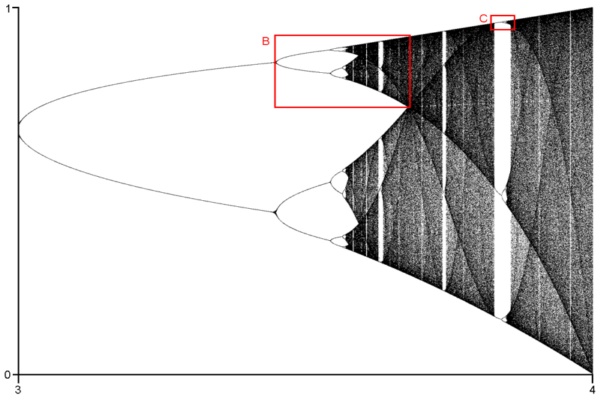}\\
\textit{Figure 1. Bifurcation diagram of the logistic map }
\end{center}

\subsection{Period-doubling bifurcations}
\hspace{.45in} The logistic map has two fixed points: $0$ and $(a-1)/a$. If $a<1$, then $|f_{a}'(0)|<1$  and we say that the fixed point $x=0$ is stable for $a<1$. If $a>1$, then $|f_a'(0)|>1$ and so the fixed point $x=0$ is now unstable and a new stable fixed point $x=(a-1)/a$ emerges. The fixed point $x=(a-1)/a$ is stable for $1<a<3$. At $a=3$, the fixed point  $x=(a-1)/a$ becomes unstable and a stable $2$-cycle (period $2$ attractor) $\{x, y\}=\{(a+1-\sqrt{a^2-2a-3})/2a,\;(a+1+\sqrt{a^2-2a-3})/2a \}$ is said to bifurcate from this fixed point. The above-mentioned $2$-cycle is the solution of the system of equations given by $y=ax(1-x)$ and $x=ay(1-y)$. When $a=1+\sqrt{6}\thickapprox 3.44949$, this $2$-cycle loses its stability and a stable $4$-cycle begins followed by a stable $8$-cycle (between $a\thickapprox 3.54409$ and $a\thickapprox 3.56440$). The period continues doubling over shorter intervals until $a\thickapprox3.56995$ where the chaotic regime takes over. A period $3$ attractor is born at $a=1+\sqrt{8}\thickapprox 3.82843$ \cite{dcMIT}, causing a notably large window in the bifurcation diagram (see Figure 1). This $3$-cycle also undergoes a period doubling cascade in which period $3\cdot2^{\ell}$ attractors are successively produced until chaos. This mechanism by which each $(k2^\ell)$-cycle loses its stability and is replaced by a stable $(k2^{\ell+1})$-cycle is called \textbf{period-doubling bifurcation} (see Figure 2).
\vspace*{-0.1in}
\begin{center}
\begin{multicols}{2}
  \includegraphics[width=2.23in]{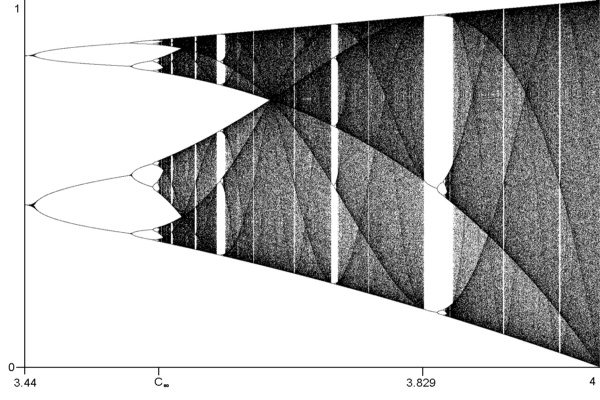}\\
  \includegraphics[width=2in]{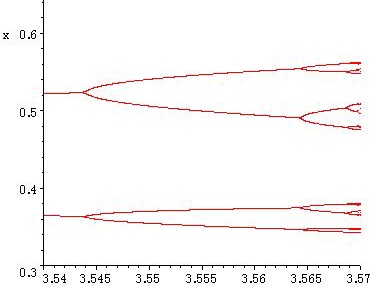}\\
 \end{multicols}
\textit{Figure 2. Magnification of the first few period doubling}
\end{center}

\indent To determine period $4$, period $8$, and in general, period $n$ solutions of the logistic map, we solve the system of difference equations, $x_{1}=x_2, x_2=x_3,\ldots ,$  and $x_n=x_1$. We note that these computations become very complicated even for relatively small values of $n$. Finding  the location of these bifurcation points on the bifurcation diagram also becomes very difficult. Even so, examining the types of orbits and their characteristic polynomials turns out to be less tortuous. This is illustrated in the following sections.

\subsection{Period-doubling bifurcation beginning with nontrivial fixed point}
$\ \ \ \ $
%\subparagraph{(12)}
The attracting period $2$ orbit of type $\theta_1=(12)$
bifurcates from the fixed point of type $\theta_0=(1)$ at $a=3$.
$\theta_1$ is the only one immediate successor to $\theta_0$, and
it follows that
$$
\begin{array}{rcl}
P_{\theta_1}(\lambda)
&=&\lambda-1\\
&=&(\lambda-1)\cdot P_{\theta_0}(\lambda),
\end{array}
$$
if we define $P_{\theta_0}(\lambda):=1$ for $\theta_0=(1)$.
Indeed the graph $G(\theta_1)=({\cal V}_1,{\cal E}_1)$ consists of
a vertex $v_{11}$ with the minus sign and
a directed edge from $v_{11}$ to itself (a loop associated with $v_{11}$):
\begin{center}
\unitlength=0.8mm
\begin{picture}(100,15)
\put(50,10){\circle*{2}}
\put(45,5){$v_{11}^{-}$}

\put(55,10){\circle{10}}
\put(60,10){\vector(0,-1){0}}
\end{picture}
\end{center}

%\subparagraph{(1324)}
The attracting period $4$ orbit of type $\theta_2=(1324)$
bifurcates from the period $2$ orbit of type $\theta_1$ at $a=1+\sqrt{6}\approx3.44949$.
$\theta_2$ with modality $+1$ is one of $2^{2-1}=2$ immediate successors to $\theta_1$.
Here we note that
$(1423)$ with modality $-1$ is another immediate successor to $\theta_1$,
however, the logistic map can realize only $+1$-modal (unimodal) cycles.
Now it follows that
$$
\begin{array}{rcl}
P_{\theta_2}(\lambda)
&=&\lambda^3-\lambda^2-\lambda+1\\
&=&(\lambda^2-1)(\lambda-1)\\
&=&(\lambda^2-1)\cdot P_{\theta_1}(\lambda).\\
\end{array}
$$
Indeed the graph $G(\theta_2)=({\cal V}_2,{\cal E}_2)$ is given as follows:
\begin{center}
\unitlength=0.8mm
\begin{picture}(100,40)
\put(45,10){\circle*{2}}
\put(65,10){\circle*{2}}
\put(55,30){\circle*{2}}

\put(42,5){$v_{23}^{-}$}
\put(55,33){$v_{21}^{+}$}
\put(60,5){$v_{22}^{-}$}

\put(63, 10){\vector(-1,0){16}}
\put(70,10){\circle{10}}
\put(75,10){\vector(0,-1){0}}
\put(53,28){\vector(-1,-2){8}}
\put(47,12){\vector(1,2){8}}
\end{picture}
\end{center}
where ${\cal V}_2=\{v_{21}^{+},v_{22}^{-},v_{23}^{-}\}$,
${\cal V}_{21}=\{v_{22}^{-}\}$, and
${\cal V}_{22}=\{v_{21}^{+},v_{23}^{-}\}$.
Here we note that
$(i)$ $G(\theta_2)|_{{\cal V}_{21}}$ is isomorphic to
$G(\theta_1)$,
the sign of vertex $v_{11}$ is consistent with that of vertex $v_{22}$,
$(ii)$ there are no directed edges from ${\cal V}_{22}$ to ${\cal V}_{21}$,
$(iii)$ $G(\theta_2)|_{{\cal V}_{22}}$ consists of
a $2$-cycle induced by $\theta_1$.

%\subparagraph{(15472638)}
The attracting period $8$ orbit of type $\theta_3=(15472638)$
bifurcates from the period $4$ orbit of type $\theta_2$ at $a\approx3.544 09$.
$\theta_3$ with modality $+1$ is
one of $2^{4-1}=8$ immediate successors to $\theta_2$.
There are other $7$ immediate successors to $\theta_2$,
however, their modalities are not equal to $+1$. Now it follows that
$$
\begin{array}{rcl}
P_{\theta_3}(\lambda)
&=&\lambda^7-\lambda^6-\lambda^5+\lambda^4-\lambda^3+\lambda+\lambda-1\\
&=&(\lambda^4-1)(\lambda^2-1)(\lambda-1)\\
&=&(\lambda^4-1)\cdot P_{\theta_2}(\lambda).\\
\end{array}
$$
Indeed the graph $G(\theta_3)=({\cal V}_3,{\cal E}_3)$ is given as follows:
\begin{center}
\unitlength=0.8mm
\begin{picture}(100,70)
\put(50,10){\circle*{2}}
\put(30,20){\circle*{2}}
\put(70,20){\circle*{2}}
\put(20,40){\circle*{2}}
\put(80,40){\circle*{2}}
\put(40,60){\circle*{2}}
\put(60,60){\circle*{2}}

\put(58,65){$v_{31}^{+}$}
\put(38,65){$v_{37}^{-}$}
\put(83,42){$v_{32}^{+}$}
\put(13,42){$v_{36}^{-}$}
\put(73,22){$v_{33}^{-}$}
\put(21,22){$v_{35}^{-}$}
\put(48,05){$v_{34}^{-}$}

\put(58,58){\vector(-3,-4){27}}
\put(78,40){\vector(-1,0){56}}
\put(78,42){\vector(-2,1){36}}
\put(68,22){\vector(-3,4){27}}
\put(50,05){\circle{10}}
\put(55,05){\vector(0,-1){0}}
\put(48,11){\vector(-2,1){16}}
\put(49,12){\vector(-1,1){27}}
\put(32,20){\vector(1,0){36}}
\put(22,41){\vector(1,0){56}}
\put(43,61){\vector(1,0){15}}
\end{picture}
\end{center}
where ${\cal V}_3=\{v_{31}^{+},v_{32}^{+},v_{33}^{-},v_{34}^{-},v_{35}^{-},v_{36}^{-},v_{37}^{-}\}$,
${\cal V}_{31}=\{v_{32}^{+},v_{34}^{-},v_{36}^{-}\}$, and
${\cal V}_{32}=\{v_{31}^{+},v_{33}^{-},v_{35}^{-},v_{37}^{-}\}$.
Here we note that
$(i)$ $G(\theta_3)|_{{\cal V}_{31}}$ is isomorphic to
$G(\theta_2)$,
the sign of vertex $v_{3(2j)}$ is consistent with that of vertex $v_{2j}$ $(j=1,2,3)$,
$(ii)$ there are no directed edges from ${\cal V}_{32}$ to ${\cal V}_{31}$,
$(iii)$ $G(\theta_3)|_{{\cal V}_{32}}$ consists of
a $4$-cycle induced by $\theta_2$.

%\subparagraph{(1,9,8,13,4,12,5,15,2,10,7,14,3,11,6,16)}
%The stable $16$-cycle of type $\theta_4=(1,9,8,13,4,12,5,15,2,10,7,14,3,11,6,16)$ bifurcates from the $8$-cycle of type $\theta_3$ at $a\approx3.56440$.$\theta_4$ with modality $+1$ is one of $2^{8-1}=128$ immediate successors to $\theta_3$. There are other $127$ immediate successors to $\theta_3$, however, their modalities are not equal to $+1$. Now it holds that
%$$\begin{array}{rcl}P_{\theta_4}(x)&=&x^{15}-x^{14}-x^{13}+x^{12}-x^{11}+x^{10}+x^9-x^8-x^7+x^6+x^5-x^4+x^3-x^2-x+1\\&=&(x^8-1)(x^4-1)(x^2-1)(x-1)\\&=&(x^8-1)\cdot P_{\theta_3}(x).\\\end{array}$$

%\subparagraph{Inductively}
Inductively,
for
the attracting period $2^{\ell}$ orbit of type $\theta_{\ell}$ with modality $+1$,
which is one of $2^{(2^{\ell-1}-1)}$ immediate successors to $\theta_{\ell-1}$,
and
bifurcated from the period $2^{\ell-1}$ orbit of type $\theta_{\ell-1}$,
we obtain $G(\theta_{\ell})$ and $P_{\theta_{\ell}}(\lambda)$.

\subsection{Period-doubling bifurcation beginning with period $3$}
$\ \ \ \ $
%\subparagraph{(123)}
Next consider the period-doubling bifurcation beginning with
period $3$ orbit of type $\theta_0=(123)$.
The graph $G(\theta_0)=({\cal V}_0,{\cal E}_0)$ consists of two vertices with a $2$-cycle and a loop:
\begin{center}
\unitlength=0.8mm
\begin{picture}(100,20)
\put(45,10){\circle*{2}}
\put(65,10){\circle*{2}}
\put(42,5){$v_{01}^{+}$}
\put(62,5){$v_{02}^{-}$}
\put(63, 9){\vector(-1,0){16}}
\put(47,11){\vector( 1,0){16}}
\put(70,10){\circle{10}}
\put(75,10){\vector(0,1){0}}
\end{picture}
\end{center}
and thus $P_{\theta_0}(\lambda)=\lambda^2-\lambda-1$.

%\subparagraph{(135246)}
The attracting period $6$ orbit of type $\theta_1=(135246)$
bifurcates from the period $3$ orbit of type $\theta_0=(123)$ at $a\approx3.8415$.
$\theta_1$ with modality $+1$ is
one of $2^{3-1}=4$ immediate successors to $\theta_0$.
$(146235),(136245),(145236)$ are also immediate successors to $\theta_0$,
however, their modalities are not equal to $+1$. Now it follows that
$$
\begin{array}{rcl}
P_{\theta_1}(\lambda)
&=&\lambda^5-\lambda^4-\lambda^3-\lambda^2+\lambda+1\\
&=&(\lambda^3-1)(\lambda^2-\lambda-1)\\
&=&(\lambda^3-1)\cdot P_{\theta_0}(\lambda)\\
\end{array}
$$
Indeed the graph $G(\theta_1)=({\cal V}_1,{\cal E}_1)$ is given as follows:
\begin{center}
\unitlength=0.8mm
\begin{picture}(100,60)
\put(50,10){\circle*{2}}
\put(20,30){\circle*{2}}
\put(80,30){\circle*{2}}
\put(40,50){\circle*{2}}
\put(60,50){\circle*{2}}

\put(58,55){$v_{11}^{+}$}
\put(38,55){$v_{15}^{-}$}
\put(83,30){$v_{12}^{+}$}
\put(13,30){$v_{14}^{-}$}
\put(48,05){$v_{13}^{+}$}

\put(59,48){\vector(-1,-4){09}}
\put(43,51){\vector(1,0){15}}
\put(78,31){\vector(-1,0){56}}
\put(22,29){\vector( 1,0){56}}
\put(23,33){\vector(1,1){16}}
\put(15,30){\circle{10}}
\put(10,30){\vector(0,1){0}}
\put(22,28){\vector(3,-2){26}}
\put(49,12){\vector(-1,4){09}}
\end{picture}
\end{center}
where
${\cal V}_1=\{v_{11}^{+},v_{12}^{+},v_{13}^{+},v_{14}^{-},v_{15}^{-}\}$,
${\cal V}_{11}=\{v_{12}^{+},v_{14}^{-}\}$, and
${\cal V}_{12}=\{v_{11}^{+},v_{13}^{+},v_{15}^{-}\}$.
Here we note that
$(i)$ $G(\theta_1)|_{{\cal V}_{11}}$ is isomorphic to
$G(\theta_0)$,
the sign of vertex $v_{1(2j)}$ is consistent with that of vertex $v_{0j}$ $(j=1,2)$,
$(ii)$ there are no directed edges from ${\cal V}_{12}$ to ${\cal V}_{11}$,
$(iii)$ $G(\theta_1)|_{{\cal V}_{12}}$ consists of
a $3$-cycle induced by $\theta_0$.

%\subparagraph{(1,5,9,4,8,11,2,6,10,3,7,12)}
%The attracting period $12$ orbit of type $\theta_2=(1,5,9,4,8,11,2,6,10,3,7,12)$ bifurcates from the period $6$ orbit of type $\theta_1$. $\theta_2$ with modality $+1$ is one of $2^{6-1}=32$ immediate successors to $\theta_1$. There are $31$ other immediate successors to $\theta_1$, however, their modalities are not equal to $+1$. Now it holds that
%$$\begin{array}{rcl}P_{\theta_2}(x)&=&x^{11}-x^{10}-x^{9}-x^{8}+x^{7}+x^{6}-x^{5}+x^{4}+x^{3}+x^{2}-x-1\\&=&(x^6-1)(x^3-1)(x^2-x-1)\\&=&(x^6-1)P_{\theta_1}(x).\\\end{array}$$
%Indeed $G(\theta_{2})$, the graph of $\theta_{2}$ is consisting of the set of vertices ${\cal X}_{2}$ and the set of paths ${\cal U}_{2}$ where ${\cal X}_{2}=\{x_{21},x_{22},x_{23},\ldots,x_{2(10)},x_{2(11)}\}$, ${\cal X}_{21}=\{x_{22},\ldots,x_{2(10)}\}$, and ${\cal X}_{22}=\{x_{21},x_{23},\ldots,x_{2(11)}\}$. Then we notice that $(i)$ $G(\theta_2)|_{{\cal X}_{21}}\cong G(\theta_1)$, the sign of vertex $x_{2(2i)}$ is consistent with that of vertex $x_{1i}$ $(i=1,2,\ldots,5)$, $(ii)$ there are no paths from ${\cal X}_{22}$ to ${\cal X}_{21}$, $(iii)$ $G(\theta_2)|_{{\cal X}_{22}}$ consists of a $6$-cycle induced by $\theta_1$.\\

%\subparagraph{Inductively}
Inductively,
for
the attracting period $(3\cdot2^{\ell})$ orbit of type $\theta_{\ell}$ with modality $+1$,
which is one of $2^{(3\cdot2^{\ell-1}-1)}$ immediate successors to $\theta_{\ell-1}$,
and
bifurcated from the period $2^{\ell-1}$ orbit of type $\theta_{\ell-1}$,
we obtain $G(\theta_{\ell})$ and $P_{\theta_{\ell}}(\lambda)$.

\subsection{Period-doubling bifurcation beginning with period $k$}
$\ \ \ \ $It is known that periodic attractors exist for all odd integers $k$.
Therefore, in general,
for the attracting period $(k\cdot2^{\ell})$ orbit of type $\theta_{\ell}$ with modality $+1$,
which is one of $2^{(k\cdot2^{\ell-1}-1)}$ immediate successors to $\theta_{\ell-1}$,
and
bifurcated from the period $2^{\ell-1}$ orbit of type $\theta_{\ell-1}$,
it follows that
$$
P_{\theta_{\ell}}(\lambda)=P_{\theta_{0}}(\lambda)\prod_{i=1}^{\ell}(\lambda^{(k\cdot2^{i-1})}-1),
$$
where $P_{\theta_{0}}(\lambda)$ is
the characteristic polynomial of
the induced matrix of $\theta_{0}$,
and $\theta_{0}$ is the type of the original period $k$ orbit.

Indeed the graph $G(\theta_{\ell})$ consists of
the set of vertices ${\cal V}_{\ell}$ and the set of directed edges ${\cal E}_{\ell}$
where
${\cal V}_{\ell}=\{v_{\ell1},v_{\ell2},v_{\ell3},\ldots,v_{\ell(k\cdot2^{\ell}-2)},v_{\ell(k\cdot2^{\ell}-1)}\}$,
${\cal V}_{\ell1}=\{v_{\ell2},\ldots,v_{\ell(k\cdot2^{\ell}-2)}\}$, and
${\cal V}_{\ell2}=\{v_{\ell1},v_{\ell3},\ldots,v_{\ell(k\cdot2^{\ell}-1)}\}$.
Here we note that $(i)$ $G(\theta_{i})|_{{\cal V}_{i1}}\cong G(\theta_{i-1})$ for $i=1,2,\ldots,\ell$,
and the signs of the vertices are preserved by each of the above isomorphisms,
$(ii)$ there are no directed edges from ${\cal V}_{i2}$ to ${\cal V}_{i1}$ for $i=1,2,\ldots,\ell$,
$(iii)$ $G(\theta_{i})|_{{\cal V}_{i2}}$ consists of a $(k\cdot2^{i-1})$-cycle induced by $\theta_{i-1}$ for $i=1,2,\ldots,\ell$, and $G(\theta_0)$ consists of $k-1$ vertices and some directed edges.

In other words,
$G(\theta_{\ell})$ essentially consists of
$k-1$ vertices with some directed edges,
a $k$-cycle of type $\theta_0$,
a $2k$-cycle of type $\theta_1$, $\cdots$,
a $(k\cdot2^{\ell-1})$-cycle of type $\theta_{\ell-1}$.
That is the reason why
$P_{\theta_{\ell}}(x)$
can be expressed as a product of
$P_{\theta_{0}}(\lambda)$, $\lambda^{k}-1$, $\lambda^{2k}-1$, $\cdots$, $\lambda^{(k\cdot2^{\ell-1})}-1$.
Thus the spectrum of $G(\theta_{\ell})$ consists of
the spectrum of $G(\theta_{0})$ and $(k\cdot2^{i-1})$ distinct roots of $\lambda^{(k\cdot2^{i-1})}-1=0$ where $i=1,2,\ldots,\ell$.

%\subparagraph{remark}
\begin{remark}\rm{\bf }
Conversely,
we may define $G(\theta_{\ell})$, the doubling of the graph $G(\theta_{i-1})$
by properties $(i)$, $(ii)$, and $(iii)$.
\end{remark}

\section{Acknowledgement}
$\ \ \ \ $The authors wish to thank Professor Evelyn L. Tan
for referring to us \cite{graph}.

%\section{Bibliography}

\end{document}